\documentclass[12pt]{amsart}
\usepackage{amsmath,amsthm,amsfonts,amssymb,amscd}
\usepackage{psfrag}
\usepackage{epsfig}
\usepackage{amsmath}
\usepackage[latin1]{inputenc}
\usepackage{amsmath,amsthm,amsfonts,amssymb,amscd}

\newcommand{\ra}{\rightarrow}

\newcommand{\diff}{\operatorname{Diff}}
\newcommand{\dw}{\diff^{1}_{\omega}(M)}
\newcommand{\difw}{\diff^{2}_{\omega}(M)}

\def\cR{\mathcal{R}}
\def\cU{\mathcal{U}}

\theoremstyle{plain}
\newtheorem{lemma}{Lemma}[section]

\newtheorem{proposition}[lemma]{Proposition}
\newtheorem{theorem}[lemma]{Theorem}
\newtheorem{maintheorem}{Theorem}

\theoremstyle{definition}

\newtheorem{remark}[lemma]{Remark}
\newtheorem{definition}[lemma]{Definition}

\title{Partial Hyperbolicity for Symplectic Diffeomorphisms}

\thanks{The authors was partially supported by  FAPESP-Proj. Tematico
03/03107-9.
A. Tahzibi would like to thank the financial support CNPq (Projeto
universal).
V. Horita was also supported by CAPES, FAPESP (02/06531-3), and PRONEX}

\author{\sc Vanderlei Horita and Ali Tahzibi}

\date{\today}

\begin{document}

\begin{abstract}
We prove that every robustly transitive and every stably ergodic
symplectic diffeomorphism on a compact manifold admits a dominated
splitting.
In fact, these diffeomorphisms are partially hyperbolic.
\end{abstract}

 \maketitle

\section{Introduction}

In the 1960's, one of the main goals in the study of dynamical
systems was to characterize the structurally stable systems and
to verify their genericity.

The Anosov-Smale theory of uniformly hyperbolic systems has
introduced a foundation to approach these subjects.
Indeed, uniform hyperbolicity proved to be the key ingredient
characterizing structurally stable systems.

We have a good description of uniformly hyperbolic dynamics.
The Smale spectral Theorem asserts that the non wandering set
of a hyperbolic diffeomorphism splits into basic pieces.
The dynamic restricted to each piece is topologically transitive.
By the topological stability of uniformly hyperbolic
systems this transitivity property persists after small
perturbations.

However, since the end of 60s, one knows that hyperbolic systems
do not constitute a dense set of dynamics, for instance see
\cite{Ne74} and \cite{AS70}.

The coexistence of infinitely many sinks and sources in a
locally residual subset ({\em Newhouse phenomenon}) shows that
uniform hyperbolicity is not a generic property.
A remarkable result of Mañé \cite{Man82} gives the following dichotomy
between uniformly hyperbolic dynamical systems and the Newhouse
phenomenon for surface diffeomorphisms.

\medskip

\noindent {\bf Theorem (Mañé \cite{Man82}).} {\em Let $S$ be a closed surface.
 Then there is a residual subset $\cR \subset \diff^1(S)$,
 $\cR = \cR_1 \cup \cR_2,$ such that
 \begin{itemize}
 \item $f \in \cR_1$ verifies Axiom A,
 \item $f \in \cR_2$ has infinitely many sinks or sources.
 \end{itemize}
}

\medskip

Recall a diffeomorphism is {\em transitive} if there exists a
point $x$ such that its forward orbit is dense. Transitive systems
do not have neither sinks nor source. We say a diffeomorphism is
$C^r$-{\em robustly transitive} if it belongs to the
$C^r$-interior of the set of transitive diffeomorphisms.

Anosov diffeomorphisms whose nonwandering set is the whole
manifold are examples of $C^1$-robustly transitive
diffeomorphisms. A direct consequence of the Mañé's Theorem is
that every $C^1-$robustly transitive diffeomorphism on a closed
surface is an Anosov diffeomorphism.

In compact manifolds of dimension greater than two, robust
transitivity does not implies uniform hyperbolicity, see
\cite{Sh71,Man78,BD96,BoV00,Ta02}. However, all these examples
present a weak form of hyperbolicity called {\em dominated
splitting}. Let $M$ be a compact manifold and $f : M \rightarrow
M$ a diffeomorphism. A $Df$-invariant decomposition $TM = E \oplus
F$ of the tangent bundle of  $M$, is dominated if for every
positive integer $\ell$ and any $x$ in $M$,
$$
\| Df^\ell|E(x) \| \cdot \| Df^{-\ell}|F(f^\ell(x)) \| < C \lambda^\ell ,
$$
for some constants $C>0$ and $0 < \lambda < 1$.

If for the dominated splitting $TM = E \oplus F$ at least one of
the subbundles is uniformly hyperbolic then $f$ is called
{\it partially hyperbolic}.

Díaz, Pujals and Ures in \cite{DPU99} show that every $C^1$-robustly
transitive diffeomorphism in a $3$-dimensional compact manifold
should be partially hyperbolic.
This result can not be generalized for higher dimension: Bonatti and
Viana present in \cite{BoV00} an example in $4$-dimensional compact
manifold of a non partially hyperbolic $C^1$-robustly transitive
diffeomorphism.
This example can be extended in any dimension, see \cite{Ta02}.

Bonatti, Díaz and Pujals in \cite{BDP03} proved that every
$C^1$-robustly transitive diffeomorphism on a $n$-dimensional
compact manifold, $n\ge 1$ admits a dominated splitting. Recently
Vivier  proved similar results for flows in \cite{V03}. She proved
that robustly transitive $C^1-$vector fields on compact manifold
do not admit singularity and that robust transitivity implies the
existence of dominated structure.

A natural question arisen for the above results is the following:
{\em What occurs if we restrict ourselves to the volume preserving
or even symplectic diffeomorphisms ?}

Let $(M,\omega)$ be a $2N$-dimensional symplectic manifold, where
$\omega$ is a non-degenerated symplectic form. We denote
$\diff^r_\omega(M)$, $r \ge 1$, the set of $C^r$-symplectic
diffeomorphisms. We say $f \in \diff^r_{\omega}(M)$ is {\em
$C^r$-robustly transitive symplectic diffeomorphism} if there
exists a neighborhood $\cU \subset \diff^r_\omega(M)$ of $f$ such
that every $g \in \cU$ is also transitive.

\begin{maintheorem}
\label{teo_transitive} Every $C^1$-robustly transitive symplectic
diffeomorphism on a $2N$-dimensional, $N\ge 1$, compact manifold
is partially hyperbolic.
\end{maintheorem}

We emphasize that if $f \in \diff_\omega^1$ and is $C^1$-robustly
transitive symplectic diffeomorphism then just symplectic nearby
diffeomorphism are transitive. So, our theorem is not a
consequence of results in \cite{BDP03}.

In the context of volume preserving  diffeomorphisms, ergodicity
of the Lebesgue measure is a basic feature. Recall that a
diffeomorphism in $\diff_\omega^1(M)$ preserves the $2$-form
$\omega$ and consequently  the volume form $\omega \wedge \dots
\wedge \omega$ is also preserved. This volume form induces in a
natural way a Lebesgue measure defined on $M$.

The stable ergodicity of symplectic diffeomorphism is defined as
follows. A symplectic diffeomorphism $f \in \diff_{\omega}^2(M)$
is {\em $C^1$-stably ergodic} if there exists a neighborhood $\cU
\subset \diff_{\omega}^1(M)$ of $f$ such that any $g \in \cU \cap
\diff_{\omega}^2(M)$ is ergodic.

The theories of stable ergodicity and robust transitivity had very
parallel development in last few years. Like as in the robust
transitivity case we know examples of stably ergodic
diffeomorphisms with weak form of hyperbolicity. Pugh and Shub
proposed to study the ergodicity among the partially hyperbolic
systems. They summarize their main theme as follows.
\begin{center}
 {\it A little hyperbolicity goes a long way in guaranteeing the
 ergodicity}
\end{center}
 We propose the reader to see
\cite{BPSW01} for approaches to prove stable ergodicity in the
partially hyperbolic case. However, there are stably ergodic
diffeomorphisms which are not partially hyperbolic, see
\cite{Ta02}.

A natural question is whether stable ergodicity of volume
preserving or even symplectic diffeomorphisms implies the
existence of a dominated splitting of the tangent bundle. For the
symplectic case we give an affirmative answer. In fact, in the
symplectic setting dominated splitting implies strong partial
hyperbolicity. This fact was first observed by Mañé see
\cite{Man84}. A proof is given in \cite{Arn02}, for $\dim M = 4$,
and in \cite{BoV03}, for the general case.

\begin{maintheorem}
\label{teo_ergodic} Every $C^1$-stably ergodic symplectic
diffeomorphism on a $2N$-dimensional, $d\ge 1$, compact manifold
is partially hyperbolic.
\end{maintheorem}

To prove our theorems we follow the arguments in \cite{BDP03}
where they obtain a dichotomy between dominated splitting and the
Newhouse phenomenon. This is done by showing that the lack of
dominated splitting leads to creation of  sinks or sources by a
convenient perturbation. Of course a symplectic diffeomorphism
does not admit sink or source. The idea in the symplectic case is
to make a perturbation and create a totally elliptic periodic
point.

In Section \ref{s.dichotomy} we state this dichotomy for symplectic
diffeomorphisms in Theorem \ref{t.dichotomy}.

We use this dichotomy to finish the proof of our main results in
Section \ref{s.proofs}. The idea is to eliminate the possibility
of creation of totally elliptic periodic points for robustly
transitive or stably ergodic symplectic diffeomorphisms. For this
purpose in Section \ref{s.proofs} we use generating functions for
symplectic diffeomorphisms to prove that ``stably ergodic and
robustly transitive symplectic diffeomorphisms are $C^1$ far from
having totally elliptic points".

To prove the dichotomy for symplectic diffeomorphisms  we prove a
similar result for linear symplectic systems. Symplectic linear
systems are introduced in Section \ref{s.linear_systems} and some
perturbation results are proved there which will be used in the rest of
the paper.

The difficulty to get the dichotomy in the symplectic case  is
that we have much less space to perform symplectic perturbations.
The perturbations have to preserve the symplectic form $\omega$.

In Section \ref{s.dichotomy_linear} we state the precise result
of dichotomy for linear symplectic systems and in Sections
\ref{s.2.4} and \ref{s.2.5} we prove this statements.
In these two final sections we prove new results for symplectic
linear systems which enable us to adapt the approach of \cite{BDP03}
for symplectic linear systems.

We mention that for robustly transitive volume preserving
diffeomorphisms the dichotomy as mentioned above is straightforward
from the result in \cite{BDP03}.
Arbieto-Matheus \cite{AM03} use this dichotomy and prove that
robustly transitive conservative diffeomorphisms have dominated
splitting.
So, the main difficulty in the proof of similar to our results in
the volume preserving case  is to show that the robustly transitive
conservative diffeomorphisms can not have totally elliptic points.
They overcome this difficulty with a new ``Pasting Lemma" which
uses a theorem  of Dacorogna-Moser \cite{DM90}.

\section{A dichotomy for symplectic diffeomorphisms}
\label{s.dichotomy}

An invariant decomposition $TM = E \oplus F$ is called
{\it $\ell$-dominated} if for any $n \ge \ell$ and every $x$
in $M$,
$$
\| Df^n|E(x) \| \cdot \| Df^{-n}|F(f^n(x)) \| < \frac{1}{2}.
$$

>From now on we use the notation $E \prec_{\ell} F$ for
$\ell$-dominated splitting and $E \prec F$ for dominated splitting
without specifying the strength of the splitting.

It is easy to verify from  definition that $\ell$-dominated
splitting has the following properties:
\begin{itemize}
\item[$1.$] If $\Lambda \subset M$ is an invariant subset that
admits an $\ell$-dominated splitting then the same is true for the
closure of $\Lambda$.
\item[$2.$] If a sequence $(f_n)_n$ of maps admitting an $\ell$-dominated
splitting converges to $f$ in $C^1$-topology then $f$ also admits an
$\ell$-dominated splitting.
\end{itemize}

\begin{theorem}
\label{t.dichotomy} Let $f\in \dw$ be a symplectic diffeomorphism
of a $2N$-dimensional manifold $M$. Then there is $\ell \in
\mathbb{N}$ such that,
\begin{itemize}
\item[a)] either there is a symplectic $\varepsilon$-$C^1$-perturbation $g$
of $f$ having a periodic point $x$ of period $n \in \mathbb{N}$ such that
$Dg^n (x) = Id$,
\item[b)] or for any symplectic diffeomorphism $g$ $\varepsilon$-$C^1$-close
to $f$ and every periodic saddle $x$ of $g$ the homoclinic class $H(x,g)$
admits an $\ell$-dominated splitting.
\end{itemize}
\end{theorem}

To prove this theorem we introduce the concept of Symplectic Linear Systems
in Sections \ref{s.linear_systems}.
Then, in Section~\ref{s.dichotomy_linear}, we reduce the proof of
Theorem~\ref{t.dichotomy} to proof a similar result for that symplectic
linear systems.

\smallskip

By the following lemma we are able to extend dominated splittings
over homoclinic classes to the whole manifold for a generic
symplectic diffeomorphism.
Before stating this result let us recall a symplectic version of
Connecting Lemma due to Xia in \cite{X96}.

\begin{theorem}[Xia] \label{t.xia}
Let $M$ be a compact $n$-dimensional manifold with a symplectic or
volume form $\omega.$ Then there is a residual subset $\cR_2
\subset \diff_{\omega}^1 (M)$ such that
 if $g \in \cR_2$ and $p \in M$ are such that $p$ is a hyperbolic periodic
 point of $g$, then $W^s(p) \cap W^u(p)$ is dense in both $W^s(p)$ and $W^u(p).$
\end{theorem}
Using this result we are able to prove the following result.

\begin{lemma}
\label{l.homoclinic}
There is a residual subset $\cR \subset \diff^1_\omega(M)$ of diffeomorphisms
$f$ such that the nontrivial homoclinic classes of hyperbolic periodic
points of $f$ are dense in $M$.
\end{lemma}

\begin{proof}
The version of this lemma for conservative diffeomorphisms
is given in \cite[Lemma 7.8]{BDP03}.
We can use the same arguments adapted for the symplectic
case by means of a result of Newhouse and connecting lemma of Xia.

First of all, we recall that for a symplectic diffeomorphism the
set of recurrent points is dense in $M$.
Indeed a symplectic diffeomorphism preserves a Lebesgue measure.
Moreover, using $C^1$-Closing Lemma of Pugh and Robinson and
the Birkhoff fixed point Theorem, Newhouse~\cite[Corollary 3.2]{Ne76}
proved that there is a $C^1$-residual subset $\cR_1 \in \diff^1_{\omega}(M)$
such that for any $g \in \cR_1$ the set of hyperbolic periodic points
of $g$ is dense in $M$.

By using Theorem~\ref{t.xia} we take $\cR=\cR_1 \cap \cR_2$.
Thus, for any $g \in \cR$ the union of non-trivial homoclinic classes
is dense.
This completes the proof of lemma.
\end{proof}

\begin{remark} \label{r.dominated}
Using the above lemma and the continuity of dominated splitting
the item (b) of Theorem~\ref{t.dichotomy} can be rewritten as
follows.
\begin{itemize}
\item[b$_1)$] the manifold $M$ is the union of finitely many invariant
(by $f$) compact sets having a dominated splitting.
\end{itemize}
We remark that the invariant compact sets mentioned above are
$\Lambda_i, i < 2N$ where $\Lambda_i$ is  closure of the union of
non-trivial homoclinic classes with an $\ell-$dominated splitting
$E \oplus F$ with dimension $i$ (i.e $\dim(E) = i$).
Of course for a transitive diffeomorphism the above item is
equivalent to have a dominated splitting on the whole manifold.
\end{remark}

\section{Proof of Theorems \ref{teo_transitive} and
\ref{teo_ergodic}} \label{s.proofs}

In this section we prove Theorems \ref{teo_transitive} and \ref{teo_ergodic}
using Theorem \ref{t.dichotomy}.

\subsection{Elliptic points vs. robust transitivity}

Let us prove that the item (a) of the dichotomy given in Theorem
\ref{t.dichotomy} does not occur for $C^1$-robustly transitive
diffeomorphisms.

\begin{lemma} \label{l.generating1}
If $f\in \dw$ has a totally elliptic periodic point $p$ of period $n$
($Df^n(p) = Id$) then there exists $g\in \dw$ and $C^1$-close to $f$
such that $g$ is not transitive.
\end{lemma}

\begin{proof}
In order to simplify our arguments, let us suppose that $p$ is a fixed
point.
We use generating functions to construct $g$ in such a way that $g$
coincide to the identity map in a small neighborhood of $p.$

Let us introduce generating function as in \cite{X96}.
We fix a local coordinated system $(x_1 , \dots , x_d, y_1 , \dots, y_d)$
in a neighborhood $U$ of $p$ such that $p$ correspond to $(0, 0)$ in this
coordinate system.
Suppose $f (x, y) = (\xi(x, y), \eta(x,y))$.
The fact $f$ is symplectic implies
$$
\sum_{i=1}^d d x_i \wedge dy_i = \sum_{i=1}^d d \xi_i \wedge d\eta_i.
$$
Moreover, we may suppose that partial derivative $\partial \eta/\partial y$ of
$\eta$ with respect to $y$ is non-singular at every point of $U$.
This enable us to solve $\eta = \eta (x, y)$ to obtain $y = y (x, \eta)$.
Hence, we can define a new system of coordinates
$(x_1, \cdots , x_d, \eta_1 , \cdots , \eta_d)$.
Let $\gamma : (x, y) \rightarrow (x, \eta (x, y))$ be the map of
change of coordinates.

Since the $1$-form
$$
\alpha:= \sum_{i=1}^d \xi_i d\eta_i + y_i dx_i
$$
is closed, there exists a $C^2$ function $ S_f (x, \eta)$, defined
on a neighborhood of $(0, \eta(0, 0))$ such that $d S_f = \alpha$.
The function $S_f$ is unique up to a constant. Moreover,
$$
\frac{ \partial S_f}{ \partial x_i} = y_i \quad \text{and} \quad
\frac{\partial S_f}{\partial \eta_i} = \xi_i.
$$

Conversely, for a real $C^2$-function $S (x , \eta)$ defined on a
neighborhood of $(0, \eta(0, 0))$ such that the second partial
derivative
$$
\frac{\partial^2 S}{\partial x \partial \eta}
$$
is non-singular in this neighborhood, we define
$$
\xi_i(x, \eta)= \frac{\partial S}{\partial \eta_i} \quad
\text{and} \quad y_i(x, \eta)= \frac{\partial S}{\partial x_i}.
$$
Then, solving $\eta$ in terms of $x, y$ we find a symplectic
diffeomorphism which maps $(x, y)$ to $(\xi, \eta)$.

Observe that in the above construction the generating function of $f$
is $C^{k+1}$ whenever $f$ is $C^k$. Moreover, given a $C^{k+1}$
generating function we obtain locally a $C^k$ symplectic diffeomorphism.

We also have that $f$ is $C^1$-close to $g$ if and only if $S_f$
is $C^2$-close to $S_g$, provided they are defined on the same
domain.

Let $\rho$ be a $C^\infty$ bump function such that
$$
\rho(z) = \left\{
          \begin{array}{ll}
          1 & \text{ if } z \in \gamma (B(\beta/2)) ,\\
          0 & \text{ if } z \notin \gamma (B(\beta))
          \end{array}
\right.
$$
where $B(r)$ is the ball of radius $r$ centered in $(0,0)$
in the $(x, y)$-coordinates.

Define a $C^2$ function given by
$$
S_g:= \rho(x, \eta) S_{id} + (1-\rho(x, \eta)) S_f.
$$
More important is that, if we take $\beta > 0$ small enough, $S_g$ is
$C^2$-close to $S_f$.
Indeed, $f$ is $C^1$-close to identity in a neighborhood of the origin
and consequently, $S_{id}$ and $S_f$ are $C^2$-close enough on
$\gamma (B(\beta))$.

Finally, it is easy to see that if $S_{id}$ and $S_{f}$ are
$\varepsilon$-close on $\gamma (B(\beta))$ in the $C^2$-topology then
$$
\rho(x, \eta) S_{id} + (1-\rho(x, \eta)) S_f
$$
is $K \varepsilon$-close to $S_f$ in the $C^2$-topology, where $K$
depends on the bump function $\rho$ and $\beta.$

Let $g \in \dw$ be the corresponding symplectic diffeomorphism
to $S_g$.
As $g$ is $C^1$-close to $f$ and locally it is equal to the identity,
we conclude that $g$ can not be transitive and this conclude the proof
of the lemma.

\end{proof}

The above lemma proves that a robustly transitive diffeomorphism
can not have a totally elliptic periodic point. So, if $f$ is
robustly transitive as in Theorem \ref{teo_transitive} then by
Remark \ref{r.dominated} we can conclude that $f$ admits a dominated
splitting.
In this way we prove Theorem \ref{teo_transitive}.

\subsection{Elliptic point vs. stable ergodicity}

In this subsection we prove that any $C^1$-stably ergodic symplectic
diffeomorphism admits a dominated splitting.
Let us recall that by $C^1$-stably ergodic diffeomorphism we mean a
symplectic $C^2$-diffeomorphism such that all symplectic
$C^2$-diffeomorphisms in a $C^1$-neighborhood are ergodic.

Let $f \in \difw$ be a stably ergodic diffeomorphism and $\mathcal{U}$
be a $C^1$-neighborhood of $f$ such that any $g \in \mathcal{U} \cap \difw$
is ergodic.
We prove that any diffeomorphism inside $\cR \cap \mathcal{U}$ admits
an $\ell$-dominated splitting where $\mathcal{R}$ is given in Lemma
\ref{l.homoclinic}.
Since the $\ell$-dominated splitting property is a closed property in
$C^1$-topology we obtain a dominated splitting for $f$.
The proof is by contradiction.

Suppose that there exists $f_1 \in \cR \cap \mathcal{U}$ close to $f$.
If $f_1$ does not admit a dominated splitting then, by Theorem~\ref{t.dichotomy},
there exists $g_1 \in \dw$ with a totally elliptic periodic point.

We claim that there exists $g \in \difw$ and $C^1$-close to $g_1$ such
that $g$ is not ergodic. This gives a contradiction with the stable
ergodicity of $f.$

Just to simplify our arguments let us suppose that $g_1$ has a totally
elliptic fixed point.
Using the same technics used in the previous subsection we construct
$g_2 \in \dw$ and $C^1$-close to $g_1$ in such a way that $g_2$ coincides
with the identity map on $B(\beta)$, for some $\beta > 0$.
Note that $g_2$ is just $C^1$ and we do not get a contradiction with
the stable ergodicity of $f$.

However, from a result of Zehnder in \cite{Ze77} there exists
$g_3 \in \difw$ and $C^1$-close to $g_2$.
In order to get a contradiction, similarly to the previous subsection,
let $\tilde \rho$ be a $C^\infty$ bump function such that
$$
\tilde{\rho}(z) = \left\{
          \begin{array}{ll}
          1 & \text{ if } z \in \gamma (B(\beta/3)) \\
          0 & \text{ if } z \notin \gamma (B(\beta/2)).
          \end{array}
\right.
$$
We define
$$
S_g (x, \eta) : = \tilde \rho (x, \eta) S_{id} + (1 - \tilde
\rho (x, \eta)) S_{g_3}.
$$
Since $g_3$ is a $C^2$-diffeomorphism we have $S_g$ is
a $C^3$-function.
Moreover, $S_g$ is $C^2$-close to $S_{g_3}$.
Therefore, $S_g$ is a generating function for a symplectic diffeomorphism
$g \in \difw$ such that $g$ is $C^1$-close to $g_3$ and
consequently $C^1$-close to $f$.

By construction, $g$ coincides with the identity map on a neighborhood
of zero. Hence $g$ can not be ergodic.
This gives a contradiction, because we have supposed that $f$ is stably
ergodic.

The proof of Theorem \ref{teo_ergodic} is complete.

\section{Symplectic Linear Systems}
\label{s.linear_systems}

In this section we introduce the key ingredient in the proof
of Theorem~\ref{t.dichotomy}.
Following the techniques of Mañé, we use symplectic linear systems
enriched with transition by Bonatti, Díaz, Pujals
(see \cite[Section 1]{BDP03}) in order to prove a dichotomy
for linear systems in Proposition~\ref{p.bdp}.
However, we stress that in our case the symplectic perturbations
we have to perform is much more difficult to realize.

Let $\Sigma$ be a topological space and $f$ a homeomorphism
defined on $\Sigma$.
Consider a locally trivial vector bundle $\mathcal{E}$ over $\Sigma$
such that $\mathcal{E}(x)$ is a symplectic vector space with
anti-symmetric non-degenerated $2$-form $\omega$.
Furthermore, we require the $\dim (\mathcal{E}(x))$ does not depend on
$x$.

We denote by $\mathcal{S}(\Sigma, f, \mathcal{E})$ the set of maps
$A : \mathcal{E} \rightarrow \mathcal{E}$ such that for every $x \in \Sigma$
the induced map $A(x, \cdot)$ is a linear symplectic isomorphism from
$\mathcal{E}(x) \rightarrow \mathcal{E}(f(x))$, that is,
$$
\omega(u, v) = \omega (A(u), A(v)).
$$
Thus, $A(x,\cdot)$ belongs to
$\mathcal{L}_\omega(\mathcal{E}(x) , \mathcal{E}(f(x)))$.

For each $x \in \mathcal{E}$ the Euclidean metric on $\mathcal{E}(x)$
and $\mathcal{E}(f(x))$ induces in a natural way a norm on
$\mathcal{L}_\omega(\mathcal{E}(x) , \mathcal{E}(f(x)))$:
$$
  |B(x, .)| = \sup \{|B(x, v)| , v \in \mathcal{E}(x), |v|= 1 \}.
$$
Furthermore, for $A \in \mathcal{S}(\Sigma, f, \mathcal{E})$ we define
$
|A| = \sup_{x \in \Sigma} |A(x,\cdot)|.
$
Note that, if $A$ belongs to $\mathcal{S}(\Sigma, f, \mathcal{E})$
its inverse $A^{-1}$ belongs to $\mathcal{S}(\Sigma, f^{-1}, \mathcal{E})$.
The {\em norm} of $A \in \mathcal{S}(\Sigma, f, \mathcal{E})$ is defined
by
$
\| A \| = \sup \{ |A|,|A^{-1}|\}.
$

Let $(\Sigma, f, \mathcal{E}, A)$ be a {\em linear symplectic system}
(or {\em linear symplectic cocycle} over $f$), that is, a $4$-tuple where
$\Sigma$ is a topological space, $f$ is a homeomorphism of $\Sigma$,
$\mathcal{E}$ is an Euclidean bundle over $\Sigma$, $A$ belongs to
$\mathcal{S}(\Sigma,f,\mathcal{E})$ and $\|A\| < \infty$.
We say that $(\Sigma, f, \mathcal{E}, A)$ is {\em periodic} if any
$p \in \Sigma$ is a periodic point of $f$.

Let $(V, \omega)$ be a symplectic vector space and $W \subseteq V$ a
vector subspace of $V$. Then the {\em symplectic complement} of $W$
is given by
$$
W^{\omega} =  \{ x \in V : \omega(x, w) = 0, \quad \text{for all}
\quad w \in
 W\}.
$$
It is easy to see that $W^{\omega}$ is also a vector space and,
by definition,
\begin{enumerate}
\item If $W \subset W^{\omega}$, then $W$ is an {\em isotropic subspace}.
\item If $W \cap W^{\omega} = 0$ , then $W$ is a {\em symplectic subspace}.
\item If $W = W^{\omega}$, then $W$ is {\em Lagrangian subspace}.
\end{enumerate}

Let $(\Sigma, f, \mathcal{E}, A)$ be a diagonalizable periodic
linear symplectic system such that $\mathcal{E}(x)= \mathbb{R}^{2N}$.
In what follows, we indicate by
$\lambda_1 (x) < \lambda_2(x) < \cdots  < \lambda_{2N}(x)$ their
eigenvalues and we denote by
$E_1(x) \prec E_2(x) \prec \cdots \prec E_{2N}(x)$ their respective
eigenspaces.

We know that if $\lambda$ is an eigenvalue of a symplectic
transformation then $\lambda^{-1}$ is also an eigenvalue. So, if
$\lambda_1 < \lambda_2 < \cdots < \lambda_{2N} $ are $2N$ distinct
eigenvalues of a symplectic transformation then $\lambda_{i^*} :=
\lambda_{2N-i+1} = \lambda_i^{-1}.$

We say that $\mathcal{B} = \{e_1, e_2, \cdots, e_{2N}\}$ is a
{\em symplectic basis} for a symplectic vector space $(V,\omega)$ if
$\mathcal{B}$ is a basis of $V$ and
$$
\omega(e_i, e_j) = \left\{
                     \begin{array}{l}
                     0 \quad \text{for} \quad j \neq i^* \\
                     1 \quad \text{for} \quad j= i^* > i.
                     \end{array}
                     \right.
$$
The fact that $\omega$ is an anti-symmetric $2$-form implies that
$\omega (e_{i^*}, e_i) = -1$ for $i^* > i.$

Let  $(\Sigma, f, \mathcal{E}, A)$ be a diagonalizable symplectic
linear system.
Suppose that  $\{e_1(p), e_2(p), \cdots , e_{2N}(p)\}$ is a
symplectic basis of $\mathcal{E}(p)$ then
$\{A(e_1(p)),A( e_2(p)), \cdots , A(e_{2N}(p))\}$ is a symplectic
basis of $\mathcal{E}(f(p))$.
For simplicity of notations, we omit
the dependence on the point of vectors $e_i$.

\begin{lemma}
Let $(\Sigma, f, \mathcal{E}, A)$ be diagonalizable periodic linear
symplectic system with distinct eigenvalues
$\lambda_1 < \cdots < \lambda_{2N}$ and
$E_1 \prec \cdots \prec E_{2N}$ be the corresponding eigenspaces.
There exists a symplectic basis $\{e_1, \cdots , e_{2N}\}$
constituted by eigenvectors.
Moreover, for $j \neq i^*$, $E_i \oplus E_{j}$ is an isotropic
subspace and $E_i \oplus E_{i^*}$ is a symplectic subspace.
\end{lemma}

\begin{proof}
Let $n = n(p)$ be the period of  $p \in \Sigma$.
Let $e_i \in E_i$ and $e_j \in E_j$ be eigenvectors of $A^n$ with
respective eigenvalues $\lambda_i$ and $\lambda_j$.
Then,
$$
\omega(A^{n}(e_i),A^{n}(e_j )) = \lambda_i \lambda_j \omega (e_i,e_j).
$$
On the other hand, since $A$ is a linear symplectic system, we have
$$
\omega(A^n(e_i),A^n(e_j)) = \omega(e_i,e_j).
$$

Thus, if $j \neq i^*$ then $\lambda_i \lambda_j \neq 1$ and
consequently $\omega (e_i, e_j) = 0$.
Moreover, as $\omega$ is non-degenerate we have $\omega(e_i, e_{i^*}) \neq 0$.
So, normalizing the vectors $e_i$, $1\le i\le N$, we can choose
a new basis constituted by eigenvectors of $A$, which we still denote
by $e_i$, such that
$\omega(e_i, e_{i^*}) = 1 $ for all $1\le i \leq N$.

Other claims in the statement are direct consequence of definitions
of symplectic and isotropic subspaces.
\end{proof}

Let $E_j \oplus E_k$ be a vector subspace of a diagonalizable
symplectic vector space. Any small symplectic perturbation of
$A|_{ E_j \oplus E_k}$ is called a {\em symplectic perturbation
along} $E_j \oplus E_k$. The next lemma asserts that any
symplectic perturbation along $E_j \oplus E_k$ can be realized as
the restriction of a symplectic perturbation of $A$. More
precisely,

\begin{lemma} (Symplectic realization) \label{l.main}
Let $(\Sigma, f, \mathcal{E}, A)$ be a diagonalizable periodic linear
system as above.
Given any $\varepsilon > 0$ and $1\le j < k \le 2N$, every
$\varepsilon$-symplectic perturbation $B$ of
$$
A |_{ E_{j} \oplus E_{k}} : E_{j} \oplus E_{k} \to E_{j} \oplus E_{k}
$$
along the orbit of $x \in \Sigma$ is the restriction of a symplectic
$\varepsilon-$perturbation $\tilde{A}$ of $A$ such that
$\tilde{A}|_{ E_i} = A|_{ E_i}$ for $i  \neq j, k, j^*, k^*.$
\end{lemma}

\begin{proof}

If $E_j \oplus E_k$ is a symplectic subspace, that is $k = j^*$, then
we define $\tilde{A}|_{E_{j} \oplus E_{j^*}} = B$, $\tilde{A}|_{E_{i}} = A$,
for $i \neq j , j^*$ and we extend it linearly.
In that way, we have
$$
\omega(\tilde{A}(e_j), \tilde{A}(e_{j^*})) = \omega (B(e_j), B(e_{j^*})) =
\omega(e_j, e_{j^*}) = 1 ,
$$
and for $i \neq j, j^*$, we get
$$
 \omega(\tilde{A}(e_i), \tilde{A}(e_j)) =
 \omega(A(e_i), \alpha e_j + \beta e_{j^*}) = 0.
$$
Moreover, for $r , s \neq j, j^*$, we have
$$
\omega(\tilde{A}(e_s), \tilde{A}(e_r)) = \omega(A(e_s), A(e_r)) =
\omega (e_s, e_r).
$$
Therefore, in this case, $\tilde{A}$ is symplectic.

\smallskip

Now, we suppose $k \neq j^*$.
An important feature we have to take account is the way we
extend $B$ to $E_{j^*} \oplus E_{k^*}$.
Once we have made it, we define $\tilde{A}$ equal to $A$
when restricted to the others subspace $E_i$,
$i \neq j, k , j^*, k^*$.
Finally, we extend linearly this operator to other vectors.

There are constants
$\alpha, \beta, \gamma, \delta \in \mathbb{R}$ such that
\begin{align*}
  B(e_j) & = \alpha e_{j} + \beta e_k \\
  B(e_k) & = \gamma e_j + \delta e_k.
\end{align*}
Let us suppose $1 \le j < k \le N$. We construct a symplectic linear
system $\tilde{A}$ in this case. In the other cases are completely
analogous, by changing conveniently the sign of the constant in
$\tilde{B}$ bellow.
We denote $\Delta = \alpha \delta - \beta \gamma$ and
we define
\begin{align*}
  \tilde{B}(e_{j^*}) & = \frac \delta \Delta e_{j^*} -
                 \frac \gamma \Delta e_{k^*} \\
  \tilde{B}(e_{k^*}) & = - \frac \beta \Delta e_{j^*} +
                 \frac \alpha \Delta e_{k^*},
\end{align*}
and extend $\tilde{B}$ linearly to $E_{j^*} \oplus E_{k^*}$.

Then, we define $\tilde{A}$ as follows:
\begin{itemize}
\item $\tilde{A}|_{E_i} = A |_{E_i}, i \notin \{j, k , j^*, k^*\}$
\item $\tilde{A}|_{E_j \oplus E_k} = B$
\item $\tilde{A}|_{E_{j^*} \oplus E_{k^*}} = \tilde{B}$,
\end{itemize}
and extend $\tilde{A}$ linearly.

Note that, if  $B$ in $E_j \oplus E_k$ is a
rotation then the perturbation $\tilde{B}$ is the same rotation
 in $E_{j*} \oplus E_{k*}$.

To verify that $\tilde{A}$ is a symplectic linear system,
it is enough to show that
$$
\omega (\tilde{A}(e_r) , \tilde{A}(e_s)) = \omega(e_r, e_s)
\text{ for any } 1 \leq r, s \leq 2N.
$$
Let us begin by the case when $r=j, s=j^*$:
\begin{align*}
\omega (\tilde{A}(e_j),\tilde{A}(e_{j^*})) & =
 \omega \left(\alpha e_j + \beta e_k,\frac{\delta}{\Delta} e_{j^*} -
 \frac{\gamma}{\Delta} e_{k^*}\right) =
 \frac{\alpha\delta}{\Delta} \omega(e_j, e_{j^*}) - \\
& \quad - \frac{\beta\gamma}{\Delta} \omega(e_{k}, e_{k^*} ) -
     \frac{\alpha \gamma}{\Delta}\omega(e_j, e_{k^*} ) +
     \frac{\beta \delta}{\Delta}\omega(e_k, e_{j^*}) \\
& = \frac{\alpha\delta}{\Delta}-\frac{\beta \gamma}{\Delta} = 1
  = \omega(e_j,e_{j^*}).
\end{align*}
Similarly, we obtain
\begin{align*}
\omega (\tilde{A}(e_{j^*}),\tilde{A}(e_j)) & = \omega(e_{j^*},e_j), \\
\omega (\tilde{A}(e_k),\tilde{A}(e_{k^*})) & = \omega(e_k,e_{k^*}), \\
\omega (\tilde{A}(e_{k^*}),\tilde{A}(e_k)) & = \omega(e_{k^*},e_k).
\end{align*}

If $r=j , s=k^* $, then
\begin{align*}
\omega (\tilde{A} (e_j), \tilde{A}(e_{k^*})) & =
   \omega\left(\alpha e_j + \beta e_k,
   -\frac{\beta}{\Delta} e_{j^*} + \frac{\alpha}{\Delta} e_{k^*}\right) =
\frac{-\alpha\beta}{\Delta} \omega(e_j, e_{j^*}) + \\
& \quad + \frac{\alpha\beta}{\Delta} \omega(e_{k}, e_{k^*} ) +
     \frac{\alpha^2}{\Delta}\omega(e_j, e_{k^*} ) -
     \frac{\beta^2}{\Delta}\omega(e_k, e_{j^*}) \\
& = 0 = \omega(e_j,e_{k^*}).
\end{align*}
Analogously, we have
$$
\omega (\tilde{A}(e_k),\tilde{A}(e_{j^*})) = 0 = \omega(e_k,e_{j^*}).
$$
Hence,
\begin{align*}
\omega (\tilde{A}(e_{j^*}),\tilde{A}(e_k)) & = 0 = \omega(e_{j^*},e_k),  \\
\omega (\tilde{A}(e_{k^*}),\tilde{A}(e_j)) & = 0 = \omega(e_{k^*},e_j).
\end{align*}

The remaining cases are direct consequence of the fact that
$\tilde{A}(e_i) = A(e_i)$ belongs to $E_i$, if $i\neq j,j^*,k,k^*$.
This completes the proof.
\end{proof}

The following lemma is used in the next sections.

\begin{lemma} \label{pertu}
Let $(\Sigma, f, \mathcal{E}, A)$ be a periodic diagonalizable
linear system.
If $\tilde E_{j}(p) \in E_i(p) \oplus E_{j}(p)$ is close to
$E_{j}(p)$ then there exists a symplectic perturbation \textsf{p}
of the identity map such that $\textsf{p} (E_i(p)) = E_i(p)$ and
\textsf{p}$( \tilde E_{j}(p)) = E_{j}(p).$
\end{lemma}

\begin{proof}
By Lemma~\ref{l.main} it is enough to define $\textsf{p}$
symplectic on $E_i(p) \oplus E_{j}(p)$ close to the identity
map.

Let us suppose $E_i(p) \oplus E_{j}(p)$ is an isotropic subspace
($ j \neq i^*$). Given $e_i \in E_i$ and $\tilde{e}_j \in
\tilde{E}_j$, we define $\textsf{p}(e_i) = e_i$ and
$\textsf{p}(\tilde{e}_j)$ the projection of $\tilde{e}_j$ over
$E_j$ and extend $\textsf{p}$ to $E_i(p) \oplus E_{j}(p)$
linearly. Since $\tilde{E}_j$ is close to $E_j$, $\textsf{p}$ is
close to the identity map. Moreover, there exists $e_j \in E_j$
such that $\tilde{e}_j = \alpha e_i + \beta e_j$, for some
constants $\alpha,\beta$. Then
$$
\omega(\textsf{p}(e_i), \textsf{p}(\tilde{e}_j))= 0 = \alpha\omega(e_i,e_i) +
\beta \omega(e_i,e_j) = \omega(e_i,\tilde{e}_j).
$$

On the other hand, let $E_i(p) \oplus E_{j}(p)$ be a symplectic subspace.
We take $e_i \in E_i , e_j \in E_j$ such that $\omega(e_i, e_j)= 1$.
Let $\tilde e_j \in \tilde E_j$ be close to $e_j$
Then, there exists constants $r$ close to zero and $s$ close to $1$
such that $\tilde e_j  = r e_i+ s e_j$.
Hence, $\{e_i,\tilde{e}_j \}$ is a basis of $E_i(p) \oplus E_{j}(p)$
and $s =\omega(e_i, \tilde e_j)$.
We define $\textsf{p}$ on this base as follows:
$$
\textsf{p}(\tilde e_{j}) = e_{j}  \quad \text{ and } \quad
\textsf{p} (e_i) = s e_i.
$$
Therefore,
$$
\omega(\textsf{p}( e_{i}),\textsf{p}(\tilde e_j)) =
\omega(\omega(e_i, \tilde e_j) e_i , e_j) =
\omega(e_i, \tilde e_j).
$$

So, in both cases $\textsf{p}$ is symplectic perturbation of the
identity map.

\end{proof}

\subsection{Symplectic Transitions}
\label{ss.transitions}

Here we recall an important notion introduced in \cite{BDP03}:
the concept of transitions.
In this work we are dealing with the systems which admit
transitions.
In order to introduce this important notion let us begin with an
example, see \cite[Section 1.4]{BDP03}.
Suppose $P$ and $Q$ are saddles of the same index linked by
transverse intersection of their invariant manifolds.
The existence of a Markov partition shows that for any fixed finite
sequence of times there is a periodic point expending alternately
the times of the sequence close to $P$ and $Q$, respectively.
Moreover, the transition time between a neighborhood of $P$ and
a neighborhood of $Q$ can be chosen bounded.
This property alow us to scatter in the whole homoclinic class of
$P$ some properties of the periodic points of this class.

Now, we introduce the concept of {\em linear systems with transitions}
as in \cite{BDP03}.
Given a set $\mathcal{A}$, a {\em word} with letters in $\mathcal{A}$
is a finite sequence of elements of $\mathcal{A}$.
The product of the word $[a]= [a_1, \dots , a_n]$ by
$[b]= [b_1, \dots , b_m]$ is the word $[a_1, \dots , a_n, b_1, \dots, b_m]$.
We say a word is {\em not a power} if $[a] \neq [b]^k$ for every word
$[b]$ and $k > 1.$

Let $(\Sigma, f, \mathcal{E}, A)$ be a periodic linear system of dimension
$2N$, that is, all $x \in \Sigma$ is a periodic point for $f$ with period
$n=n(x)$. We denote $M_A$ the product $A^n(x)$ of $A$ along the orbit of
$x$.

If $(\Sigma, f, A)$ is a periodic symplectic linear system of matrices in
$SP(2N, \mathbb{R})$, then for any $x \in \Sigma$ we write,
$$
[M]_A (x) = (A(f^{n-1}(x)), \dots , A(x)),
$$
where $n$ is period of $x$.
The matrix $M_A(x)$ is the product of the words $[M]_A (x)$.

\begin{definition}[Definition 1.6 of \cite{BDP03}]
\label{d.transition} Given $\varepsilon > 0$, a periodic linear
system $(\Sigma, f, \mathcal{E},A)$ {\em admits
$\varepsilon$-transitions} if for every finite family of points
$x_1, \dots , x_n = x_1 \in \Sigma$ there is an orthonormal system
of coordinates of the linear bundle $\mathcal{E}$ so that
$(\Sigma, f, \mathcal{E},A)$ can now be considered as a system of
matrices $(\Sigma, f, A)$), and for any $(i, j) \in \{1,\dots
,n\}^2$ there exist $k(i, j) \in \mathbb{N}$ and a finite word
$[t^{i,j}]= (t_1^{i,j}, \dots , t_{k(i,j)}^{i,j})$ of matrices in
$SP(2N, \mathbb{R})$, satisfying the following properties:
\begin{enumerate}
\item For every $m \in \mathbb{N}$, $\imath = (i_1, \dots, i_m) \in \{1,\dots,n\}^m$,
and $\alpha = (\alpha_1,\dots ,\alpha_m) \in \mathbb{N}^m$
consider the word
\begin{align*}
[W(\imath, \alpha)] & = [t^{i_1,i_m}][M_A(x_{i_m})]^{\alpha_m}
[t^{i_m,i_{m-1}}][M_A(x_{i_{m-1}})]^{\alpha_{m-1}} \dots \\
& \quad \dots [t^{i_2,i_1}][M_A(x_{i_1})]^{\alpha_1},
\end{align*}
where the word
$w(\imath, \alpha) = ((x_{i_1}, \alpha_1) ,\dots,(x_{i_m}, \alpha_m))$
with letters in $M \times \mathbb{N}$ is not a power.
Then there is $x(\imath, \alpha) \in \Sigma$ such that
\begin{itemize}
\item The length of $[W(\imath, \alpha)]$ is the period of
$x(\imath,\alpha)$.
\item The word $[M]_A (x(\imath, \alpha))$ is $\varepsilon$-close to
$[W(\imath, \alpha)]$ and there is an $\varepsilon$-symplectic perturbation
$\tilde A$ of $A$ such that the word $[M]_{\tilde{A}} (x(\imath,\alpha))$
is $[W(\imath, \alpha)]$.
\end{itemize}

\item One can choose $x(\imath, \alpha)$ such that the distance between
the orbit of $x(\imath, \alpha)$ and any point $x_{i_k}$ is bounded by some
function of $\alpha_k$ which tends to zero as $\alpha_k$ goes to infinity.
\end{enumerate}
\end{definition}

Given $\imath ,\alpha$ as above, the word $[t^{i, j}]$ is an
$\varepsilon$-transition from $x_j$ to $x_i$.
We call {\em $\varepsilon$-transition matrices} the matrices $T_{i,j}$
which are product of the letters composing $[t^{i, j}]$.
We say a periodic linear system system {\em admits transitions} if
for any $\varepsilon > 0$ it admits $\varepsilon$-transitions.

The following lemma gives an example of linear systems with symplectic
transitions.
It is a symplectic version of \cite[Lemma 1.9]{BDP03} and its proof,
based on the existence of Markov Partitions, is analogous the proof of
\cite{BDP03}.

\begin{lemma} \label{l.homoclinicclasses}
Let $f$ be a symplectic diffeomorphism and let $P$ be a hyperbolic
saddle of index $k$ (dimension of its stable manifold).
The derivative $Df$ induces a continuous periodic symplectic
linear system on the set $\Sigma$ of hyperbolic saddles in the
homoclinic class $H(P,f)$ of index $k$ and homoclinically related
to $P.$
\end{lemma}

\begin{remark}
We have some good properties that we use during the next sections.
Consider points $x_1, \dots , x_n=x_1 \in \Sigma$ and
$\varepsilon$-transitions $[t^{i,j}]$ from $x_j$ to $x_i$.
Then
\begin{enumerate}
\item for every positive $\alpha \geq 0$ and $\beta \geq 0$ the word
$$
([M]_A(x_i))^{\alpha}) [t^{i, j}] ([M]_A(x_j))^{\beta})
$$
is also an $\varepsilon$-transition from $x_j$to $x_i$.

\item for any $i, j$ and $k$ the word $[t^{i, j}][t^{j, k}]$
is an $\varepsilon$-transition from $x_k$ to $x_i.$
\end{enumerate}
\end{remark}

The following lemma whose proof is analogous of \cite[Lemma 1.10]{BDP03}
states that every periodic symplectic system with transitions can be
approximated by a diagonalizable systems defined on a dense
subset of $\Sigma$.
We emphasize that this lemma is also true for symplectic case.

\begin{lemma}
\label{l.diagonalizable}
Let $(\Sigma, f, \mathcal{E}, A)$ be a periodic linear symplectic
system with transition.
Then for any $\varepsilon > 0$ there is a
diagonalizable symplectic $\varepsilon$-perturbation $\tilde{A}$ of
$A$ defined on a dense invariant subset $\tilde{\Sigma}$ of
$\Sigma.$
\end{lemma}

\begin{remark}
We remark that the diagonalizable system near to $A$
as required in the above lemma is not necessarily continuous, but
it does not matter in the way we apply this lemma.
\end{remark}

\section{A dichotomy for Symplectic linear systems}
\label{s.dichotomy_linear}

In this section, we reduce the study of the dynamics of symplectic
diffeomorphisms in Theorem~\ref{t.dichotomy} to a problem on
symplectic linear systems in Proposition~\ref{p.bdp}. We split the
proof of this proposition in two propositions (Propositions~\ref{p.bdp2.4}
and \ref{p.bdp2.5}) whose proofs is given in Sections~\ref{s.2.4}
and \ref{s.2.5}, see \cite[Section 2.1]{BDP03} for more
details.

An important tool to make the interplay between a dichotomy for
diffeomorphisms and for linear symplectic systems is a symplectic
version of Frank´s Lemma.

\begin{lemma}[Symplectic Franks' Lemma]
\label{l.franks}
Let $f \in \dw$ and $E$ a finite $f$-invariant set.
Assume that $B$ is a small symplectic perturbation of $Df$
along $E$.
Then for every neighborhood $V$ of $E$ there is a symplectic
diffeomorphism $g$ arbitrarily $C^1$-close to $f$ coinciding with
$f$ on $E$ and out of $V$, and such that $Dg$ is equal to $B$ on
$E$.
\end{lemma}

\begin{remark}
\label{r.franksk}
If the symplectic diffeomorphism in the Lemma~\ref{l.franks}
belongs to $\diff^{r}_{\omega}(M)$, $1 \le r \le \infty$ then
the corresponding perturbed diffeomorphism can be taken in
$\diff^{r}_{\omega}(M)$.
That is because the perturbation envolves exponential function and
some bump functions which are $C^\infty$.
However, it is important to note that the perturbed diffeomorphism
is just $C^1$-close to the initial one.
\end{remark}

The following proposition is a result that provides a dichotomy
for symplectic linear systems.

\begin{proposition}[Main Proposition]
\label{p.bdp} For any $\varepsilon > 0$, $N\in \mathbb{N}$ and $K
> 0$ there is $\ell>0$ such that any continuous periodic
$2N$-dimensional linear system $(\Sigma,f,\mathcal{E},A)$ bounded
by $K$ (i.e. $\|A\| < K$) and having symplectic transitions
satisfies the following,
\begin{itemize}
\item either $A$ admits an $\ell$-dominated splitting,
\item or there are a symplectic $\varepsilon$-perturbation
$\tilde{A}$ of $A$ and a point $x \in \Sigma$ such that
$M_{\tilde{A}}(X)$ is the identity matrix.
\end{itemize}
\end{proposition}

In the proof of this theorem we use the following notions.

\begin{definition}[Definition 2.2 of \cite{BDP03}]
Let $M\in GL(N,\mathbb{R})$ be a linear isomorphism of $\mathbb{R}^N$ such
that $M$ has some complex eigenvalue $\lambda$, i.e.,
$\lambda \in \mathbb{C} \setminus \mathbb{R}$. We say $\lambda$
has {\em rank} $(i,i+1)$ if there is a $M$-invariant splitting of
$\mathbb{R}^N$, $F\oplus G \oplus H$, such that:
\begin{itemize}
\item every eigenvalue $\sigma$ of $M|_F$ (resp. $M|_H$) has modulus
$|\sigma| < |\lambda|$ (resp. $|\sigma| > |\lambda|$),

\item $\dim (F) = i-1$ and $\dim(H) = N -i -1$,

\item the plane $G$ is the eigenspace of $\lambda$.
\end{itemize}
\end{definition}

We say a periodic linear system $(\Sigma, f, \mathcal{E},A)$
{\em has a complex eingenvalue of rank $(i,i+1)$} if there is
$x\in \Sigma$ such that the matrix $M_A(x)$ has a complex eigenvalue
of rank $(i,i+1)$.

\smallskip

We split the proof of Proposition~\ref{p.bdp} into the following
two results, which are symplectic version of
\cite[Propositions 2.4 and 2.5]{BDP03}:

\begin{proposition}
\label{p.bdp2.4} For every $\varepsilon > 0$, $N\in \mathbb{N}$
and $K > 0$ there is $\ell \in \mathbb{N}$ satisfying the
following property: Let $(\Sigma,f,\mathcal{E},A)$ be a continuous
periodic $2N$-dimensional linear system with symplectic
transitions such that its norm $\|A\|$ is bounded by $K$. Assume
that there exists $i \in \{1,\dots , 2N-1\}$ such that every
symplectic $\varepsilon$-perturbation $\tilde{A}$ of $A$ has no
complex eigenvalues of rank $(i,i+1)$. Then
$(\Sigma,f,\mathcal{E},A)$ admits an $\ell$-dominated splitting
$F\oplus G$, $F\prec_\ell G$, with $\dim(F) = i$.
\end{proposition}

\begin{proposition}
\label{p.bdp2.5} Let $(\Sigma,f,\mathcal{E},A)$ be a periodic
linear system with symplectic transitions. Given $\varepsilon >
\varepsilon_0 > 0$ assume that, for any $i \in \{1,\dots ,
2N-1\}$, there is a symplectic $\varepsilon_0$-perturbation of $A$
having a complex eigenvalue of rank $(i,i+1)$. Then there are a
symplectic $\varepsilon$-perturbation $\tilde{A}$ of $A$ and $x\in
\Sigma$ such that $M_{\tilde{A}}$ is the identity matrix.
\end{proposition}

The proof of Proposition \ref{p.bdp2.4} follows from $2$-dimensional
arguments of Mañé \cite{Man82} and higher dimensional arguments of
Bonatti-Díaz-Pujals in \cite{BDP03}.
However, when we change a symplectic linear system along a subspace
by a symplectic perturbation, it produces an effect along it conjugated
symplectic subspace.
So, in many time, we have to deal with $4$-dimensional arguments
instead of Mañé's $2$-dimensional arguments.

Although our arguments are based on \cite{BDP03}, all
perturbations we have to perform are symplectic. For this reason,
we have to be more careful and introduce new techniques.

\section{Proof of Proposition \ref{p.bdp2.4}} \label{s.2.4}

The main ideas of the proof of Proposition \ref{p.bdp2.4} is
to use the argument of Mañé in $2$-dimensional case and
reduction techniques in \cite{BDP03}.
In fact, by symplectic nature of our systems, the problem will
not be reduced to a $2$-dimensional problem.
Roughly speaking, the problem will be reduced to a problem
in a $4$-dimensional subspace.

So, first of all, let us recall the $2$-dimensional version of
Proposition~\ref{p.bdp2.4}.

\begin{proposition}[Mañé]
\label{p.asli}
Given any $K$ and $\varepsilon > 0$ there is $\ell \in \mathbb{N}$
such that for every $2$-dimensional linear system
$(\Sigma, f, \mathcal{E}, A)$, with norm $\|A\|$ bounded by $K$ and
such that the matrices $M_A(x)$ preserve the orientation,
\begin{enumerate}
\item either $A$ admits an $\ell$-dominated splitting,

\item or there are an $\varepsilon$-perturbation $B$ of $A$ and
$x \in \Sigma$ such that $M_{B}(x)$ has a complex (non-real)
eigenvalue.
\end{enumerate}
\end{proposition}

\begin{remark}
If $A$ in the above proposition preserves a non-degenerate form
$\omega$ defined on $\mathcal{E}$ then it is possible to choose
$B$ $\varepsilon$-close to $A$ also preserving $\omega$.
Recall that, in $2$-dimensional setting, if $\det(B) = 1$ then
$B$ preserves $\omega$.
In order to construct such $B$ suppose that $\det(\tilde{B}) \neq 1$
where $\tilde{B}$ comes from the above proposition, we just
substitute $\tilde{B}$ by $B = \tilde{B}/\det(\tilde{B})$.
The fact that the determinant map is continuous and $\det(A) = 1$
implies that $B$ is close to $A.$
\end{remark}

\subsection{Dimension Reduction}
To generalize Proposition \ref{p.asli} to higher dimensions,
Bonatti-Diaz-Pujals in \cite{BDP03} introduced dominated splittings
for quotient space.
Here we just recall the definitions, for more details we recommend
the reader to see \cite[Section 4]{BDP03}.

Let $(\Sigma, f, \mathcal{E}, A)$ be a linear system and $F$ an
invariant subbundle of $\mathcal{E}$ with constant dimension.
We denote by $A_F$ the restriction of $A$ to $F$ and by $A/F$ the
quotient of $A$ along $F$ endowed with the metric of orthogonal
complement $F^{\perp}$ of $F$; i.e., given a class $[v]$ we let
$|[v]| = |v_F^{\perp}|$.

An important statement proved in \cite{BDP03} is the following.

\begin{lemma} [Lemma 4.4 of \cite{BDP03}]
For any $K > 0$ and $\ell \in \mathbb{N}$, there exists $L$ with
the following property:
Given any linear system $(\Sigma, f, \mathcal{E}, A)$ such that
$\|A\|$ is bounded by $K$ with an invariant splitting
$E \oplus F \oplus G$, one has
\begin{enumerate}
\item $ E \prec_{\ell} F$ and $E/F \prec_{\ell} G/F \Rightarrow
E \prec_{L} F \oplus G$.

\item $F \prec_{\ell} G$ and $E/F \prec_{\ell} G/F \Rightarrow
E \oplus F \prec_{L}  G$.
\end{enumerate}
\end{lemma}

Using the previous lemma, the proof of Proposition~\ref{p.bdp2.4}
follows from Lemma~\ref{sub2.4} besides a inductive process completely
analogous to that in \cite[Lemmas 5.2 and 5.3]{BDP03}.
In fact the next lemma is a symplectic version of \cite[Lemma 5.1]{BDP03}
and it can be understood as a $2$-dimensional version of
Proposition~\ref{p.bdp2.4} and in it proof we have to take account
that we are deal with symplectic systems.
Therefore, we give a complete proof of this statement.
We omit the inductive argument necessary to the proof of
Proposition~\ref{p.bdp2.4}, since it is similar to the mentioned work.

\begin{lemma} \label{sub2.4}
Given $K>0$ and $\varepsilon > 0$ there is $\ell \in \mathbb{N}$
such that for any diagonalizable linear periodic system
$(\Sigma, f, \mathcal{E}, A)$ of dimension $2N$ and bounded by $K$,
and any $1\le i \le 2N-1$ one has
\begin{itemize}
\item Either there is an $\varepsilon$-perturbation of $A$ having
a complex eigenvalue of rank $(i,i+1)$, \item or for every $j \le
i \le k$
$$
E_j/(E_{j+1}\oplus \dots \oplus E_k) \prec_\ell E_{k+1}/(E_{j+1}
\oplus \dots \oplus E_k).
$$
\end{itemize}
\end{lemma}

\begin{proof}
Suppose that
$\lambda_1 < \dots < \lambda_N < \lambda_{N+1} < \dots < \lambda_{2N}$
are the eigenvalues of $A$.
As $A$ is symplectic we have
$\lambda_i =  (\lambda_{2N - i + 1})^{-1} = \lambda_{i^*}^{-1}. $

Fix $\varepsilon > 0$ and let $\ell$ be the dominance constant in the
Proposition~\ref{p.asli}.
If $E_j/(E_{j+1}\oplus \dots \oplus E_k) \prec_\ell E_{k+1}/(E_{j+1}
\oplus \dots \oplus E_k)$ we are done.
Otherwise, by Proposition~\ref{p.asli}, we perturb the quotient to
collapse the eigenvalues $\lambda_j$ and $\lambda_{k+1}$.
Moreover, by Lemma~\ref{l.main}, this perturbation of the quotient
gives a perturbation $\tilde{A}$ of $A$ having a pair of eigenvalues
$\tilde{\lambda}_j = \tilde{\lambda}_{k+1}$, which are respectively
continuation of $\lambda_j$ and $\lambda_{k+1}$ and preserving the
eigenvalues of the restriction of $A$ to
$E_{j+1} \oplus \dots \oplus E_k$.

Consider a symplectic isotopy $A_t$\,, $0 \le t \le 1$, such that
$A_0 = A$, $A_1 = \tilde{A}$ and denote $\lambda_{j,t}$ and
$\lambda_{k+1,t}$ the continuations of $\lambda_j$ and $\lambda_{k+1}$
at time $t$.
Let us assume that $\lambda_{j,t} \le \lambda_{k+1,t}$ for every
$0 \le t < 1$.
We analyze the following cases:
\begin{enumerate}
\item $\lambda_i \le 1 \le \lambda_{i+1}$,
\item $\lambda_i < \lambda_{i+1} < 1$,
\item $1 < \lambda_i < \lambda_{i+1}$.
\end{enumerate}

Observe that, as this perturbation is symplectic, the eigenvalues
$\tilde{\lambda}_{j^*}$ and $\tilde{\lambda}_{(k+1)^*}$ also will
collapse.
So, we get $\tilde{A}$ having the same eigenvalues $\lambda_s$,
for $s \notin \{ j, j^*, k+1, (k+1)^*\}$ of $A$ and
$\tilde{\lambda}_j = \tilde{\lambda}_{k+1}$,
$\tilde{\lambda}_{(k+1)^*} = \tilde{\lambda}_{j^*}$.
Furthermore, when we get a complex eigenvalue of rank $(i,i+1)$,
we also get a complex eigenvalue of rank $((i+1)^*,i^*)$.
Hence, the proof of item $(2)$ yields to the proof of item $(3)$.

In order to proof item $(1)$, note that we must have $i = N$.
It is because $\lambda_i \le 1$ and $\lambda_{i+1} = \lambda_i^* \ge 1$
are consecutive eigenvalues of a symplectic system.
Moreover, since $A$ is diagonalizable then $\lambda_N$ and
$\lambda_{N+1} = \lambda_N^*$ can not assume the value $1$.
For the proof of this item, we have the following alternatives:
\begin{itemize}
\item[$(1.a)$] $\tilde{\lambda}_j = \tilde{\lambda}_{k+1} < \lambda_N$.
So, there exists $0 \leq t \leq 1$ such that $\lambda_{k+1,t} = \lambda_N$
and $\lambda_{j,t}<\lambda_N$.
Hence, $\lambda_{(k+1)^*,t} = \lambda_{N+ 1} > 1$.
Therefore, there exists $t' < t$ such that
$\lambda_{k+1, t'} = \lambda_{(k+1)^*, t'} = 1$.
Then, we perturb $A_{t'}$ to get a complex eigenvalue of rank $(i, i+1)$.

\item[$(1.b)$] $\tilde{\lambda}_j = \tilde{\lambda}_{k+1} > \lambda_{N+1}$.
This case is similar to the previous one.

\item[$(1.c)$] $ \lambda_N < \tilde{\lambda}_j = \tilde{\lambda}_{k+1}
< \lambda_{N+1}$.
Recall that $\lambda_N < 1 < \lambda_{N+1}$.
So, there exists $t$ such that either $\lambda_{j, t} = \lambda_{j^*,t}=1$
or $\lambda_{k+1, t} = \lambda_{(k+1)^* , t} = 1$.
In both cases we obtain a complex eigenvalue of rank $(i, i+1)$ after a
small perturbation of $A_t$.
\end{itemize}

Now, we consider the second case where
$\lambda_i < \lambda_{i+1} < 1$.
Again we take account the following alternatives:
\begin{itemize}
\item[$(2.a)$] $\tilde{\lambda}_j = \tilde{\lambda}_{k+1} < \lambda_i$.
If $\lambda_{k + 1} < \lambda_{(i+1)^*}$ then there exists $0 \le t \le 1$
such that $\lambda_{k+1, t} = \lambda_{i}$.
After a small perturbation of $A_t$, we get a complex eigenvalue of rank
$(i,i+1)$.
Otherwise, $\lambda_{k + 1} > \lambda_{(i+1)^*}$ implies
$\lambda_{(k+1)^*} < \lambda_i$.
Then, there exists $0 \le t' < 1$ such that either
$\lambda_{j,t'} < \lambda_{(k+1)^*, t'} = \lambda_{i+1}$ (recall that
$\tilde{\lambda}_{(k+1)^*} > 1$) or
$\lambda_{(k+1)^*, t'} < \lambda_{j,t'} = \lambda_{i+1}$ (there are no
reason to $\lambda_{j,t}$ remains less than $\lambda_{i+1}$ during all
the isotopy).
In both cases we perturb slightly $A_{t'}$ to produce a complex eigenvalue
of rank $(i,i+1)$.

\item[$(2.b)$] $\tilde{\lambda}_j = \tilde{\lambda}_{k+1} > \lambda_{i+1}$.
If $\lambda_{k + 1} < \lambda_{(i+1)^*}$, let $0 \le t' \le 1$ be the smallest
$t$ such that
or $\lambda_i < \lambda_{j,t} = \lambda_{(k+1)*, t} < \lambda_{i+1}$ or
$\lambda_{j,t} = \lambda_{i+1} < \lambda_{(k+1)*,t}$ or
$\lambda_{j,t} < \lambda_{i} = \lambda_{(k+1)*,t}$.
After a small perturbation of $A_{t'}$, we get a complex eigenvalue of
rank $(i,i+1)$.
Otherwise, $\lambda_{k + 1} > \lambda_{(i+1)^*}$ implies
$\lambda_{(k+1)^*} < \lambda_{i}$.
Then, there exists $0 \le t' \le 1$ such that either
$\lambda_{(k+1)^*, t'} < \lambda_{j,t'} = \lambda_{i+1}$ or
$\lambda_{j,t'} < \lambda_{(k+1)^*, t'} = \lambda_{i+1}$.
In both cases we perturb slightly $A_{t'}$ to produce a complex
eigenvalue of rank $(i,i+1)$.

\item[$(2.c)$] $\lambda_i < \tilde{\lambda}_j = \tilde{\lambda}_{k+1} <
\lambda_{i+1}$.
If $\lambda_{k + 1} > \lambda_{(i+1)^*}$ then $\lambda_{k + 1}^* < \lambda_{i}$.
Then, there exists $0 \le t' \le 1$ such that either
$\lambda_{j,t'} < \lambda_{(k+1)^*, t'} = \lambda_i$ or
$\lambda_{(k+1)^*, t'} < \lambda_{j,t'} = \lambda_i$.
In both cases we perturb slightly $A_{t'}$ to produce a complex
eigenvalue of rank $(i,i+1)$.
Otherwise, $\lambda_{k + 1} < \lambda_{(i+1)^*}$ implies that
a small perturbation of $\tilde{A}$ gives us a complex
eigenvalue of rank $(i,i+1)$.
\end{itemize}

This completes the proof.
\end{proof}

\medskip

\noindent{\bf End of the proof of Proposition \ref{p.bdp2.4}}:
Observe that Lemma~\ref{sub2.4} is for diagonalizable systems.
But, now we use Lemma \ref{l.diagonalizable} to end the proof
of Proposition \ref{p.bdp2.4}.

Assume that there are $\varepsilon > 0$ and
$i \in \{ 1, \dots , 2N-1\}$ such that every $\varepsilon$-perturbation
of $A$ has no complex eigenvalue of rank $(i, i+1).$
Choose a sequence $\varepsilon _n < \varepsilon/2$ converging
to zero.
As the system $(\Sigma, f, A)$ has transition, using Lemma
\ref{l.diagonalizable}, we get a dense subset $\Sigma_n$ and
diagonalizable $\varepsilon_n$-perturbation $B_n$ of $A$ defined
on $\Sigma_n$.
Then, we apply Lemma~\ref{sub2.4} for $B_n$: there is an integer
$L > 0$ such that every $B_n$ admits an $L-$dominated splitting
$E_n \oplus F_n$ with $\dim(E_n)= i$.
Finally, as $\Sigma_n$ are dense and $\|B_n - A\| \ra 0$, we
conclude that $A$ admits an $L$-dominated splitting $E \oplus F$
with $\dim (E) = i$.

\section{Proof of Proposition \ref{p.bdp2.5}}
\label{s.2.5}

After getting periodic points with complex eigenvalues with rank
$(i, i+1)$ the idea in the proof of Proposition \ref{p.bdp2.5}
is to use the symplectic transitions to multiply matrices
corresponding to different points of $\Sigma$ having complex
eigenvalues of different ranks.
Observe that as $f$ is symplectic a periodic point with complex
eigenvalue of rank $(i, i+1)$ is also of rank $((i+1)^*,i^*)$.

The next proposition is a symplectic version of Lemma 5.4 of
\cite{BDP03}.

\begin{proposition}
\label{p.identity}
Let $(\Sigma, f, \mathcal{E}, A)$ be a continuous periodic symplectic
linear system with symplectic transitions.
Fix $\varepsilon_0 > 0$ and assume that a symplectic
$\varepsilon_0$-perturbation of $A$ has a complex eigenvalue of
rank $(i, i+1)$ for some $i \in \{1, \cdots , 2N-1\}$.
Then for every $0 < \varepsilon_1 < \varepsilon_0$ there is a point
$p \in \Sigma$ such that for every $1 \leq i < 2N$ there is a
symplectic $\varepsilon_1$-transition $[t^i]$ from $p$ to itself
with the following properties:
\begin{itemize}
\item There exists a symplectic $\varepsilon_1$-perturbation
$[M]_{\tilde{A}}(p)$ of the word $[M]_A(p)$ such that the
corresponding matrix $[M]_{\tilde{A}}(p)$ has only real positive
eigenvalues with multiplicity $1$.
Denote by $\tilde{\lambda}_1 < \dots < \tilde{\lambda}_{2N}$ such
eigenvalues and by $E_i(p)$ their respective ($1$-dimensional)
eigenspaces.

\item There is a symplectic $(\varepsilon_0+\varepsilon_1)$-perturbation
$[\tilde{t}^i]$ of the transition $[t^i]$ such that the corresponding
matrix $\tilde{T}^i$  satisfies
\begin{itemize}
\item $\tilde{T}^i(E_j(p)) = E_j(p) \text{ if } j
  \notin \{i,i+1, i^*, (i+1)^*\}$,
\item $\tilde{T}^i(E_i(p)) = E_{i+1}(p) \text{ and }
\tilde{T}^i(E_{i+1}(p)) = E_i(p)$.
\end{itemize}
\end{itemize}
\end{proposition}

A key tool in the proof of Proposition~\ref{p.identity} is
symplectic transitions constructed in Proposition~\ref{p.ultima}.
These transitions preserve the dominated splitting corresponding
to two different periodic points.

Let $p, p_i \in \Sigma$ such that
$$
\mathcal{E}(p)= F_1 \prec \dots \prec F_{2N},
$$
and
$$
\mathcal{E}(p_i) = E_1 \prec \dots \prec E_{(i+1)^*, i^*} \prec
\dots \prec E_{i, i+1} \prec \dots \prec E_{2N},
$$
where $F_i$ and $E_i$ are $1$-dimensional eigenspaces and
$E_{i, i+1}$ and $E_{(i+1)^*, i^*}$ are $2$-dimensional eigenspaces
corresponding to the complex eigenvalues.
We fix a symplectic basis $\{f_1, \dots , f_{2N}\}$ for
$\mathcal{E}(p)$.
The eigenvalues corresponding to $E_{(i+1)^*, i^*}$ are inverse
of the eigenvalues corresponding to $E_{i, i+1}$.
We only deal with the case $i > N$, that is when $E_{i, i+1}$ appears
after $E_{(i+1)^*, i^*}$ in the above dominated splitting.
The another case is completely similar to this case.

\begin{proposition}
\label{p.ultima}
There exists a symplectic transition $[t_{i,0}]$ from $p$ to $p_i$ such
that its matrix $T_{i,0}$ satisfies the following:
\begin{itemize}
\item $T_{i, 0}(F_j) = E_j$, if $j \neq i, i+1, i^*,(i+1)^*$,
\item $T_{i, 0}(F_i \oplus F_{i+1}) = E_{i,i+1}$, and
\item $T_{i, 0}(F_{(i+1)^*} \oplus F_{i^*}) = E_{(i+1)^*,i^*}$.
\end{itemize}
There is also a symplectic transition $[t_{0,i}]$ from $p_i$ to $p$ with a similar
properties.
\end{proposition}

To prove the proposition we state two auxiliary lemmas.

\begin{lemma}
\label{l.1}
Under the hypotheses of the proposition \ref{p.ultima}, if $i=2N-1$
there exists a symplectic transition  from $p$ to $p_i$ such that its matrix
$T$ satisfies the following properties:
\begin{itemize}
\item  $T(F_{2N-1} \oplus F_{2N}) = E_{2N-1,2N}$ and
$T (F_{1} \oplus F_{2}) = E_{1,2},$

\item $T(F_3 \oplus \cdots \oplus F_{2N-2}) = E_3 \oplus \dots
\oplus E_{2N-2}$.
\end{itemize}

\end{lemma}

\begin{proof}
Suppose that $t_{i,0}$ is an arbitrary transition from $p$ to $p_i$.
Let $M_i$ and $M$ denote respectively $M_A(p_i)$ and $M_A(p).$
After a small symplectic perturbation we may suppose
$$
T_{i, 0} (F_{2N-1}\oplus F_{2N}) \nsubseteq E_{1,2} \oplus E_3
\cdots \oplus E_{2N-2}.
$$

\begin{figure}[htbp]
  \centering
  \psfrag{Ti,0}{$T_{i,0} (F_{2N-2} \oplus F_{2N})$}
  \psfrag{T0i}{$T_{0,i}$}
  \psfrag{F2N-1}{$F_{2N-1}$}
   \psfrag{F2N-1,2N}{$E_{2N-1,2N}$}
   \psfrag{F2N}{$F_{2N}$}
 \psfrag{p}{$p$}
  \psfrag{Min}{$M_i^{n_1} (T_{i,0} (F_{2N-1} \oplus F_{2N}))$}
  \psfrag{pi}{$p_i$}
\psfrag{ pi}{$p_i$}
  \psfrag{S}{$S$}
  \psfrag{E12}{$E_{1,2}$}
  \includegraphics[width=12cm, height=6cm]{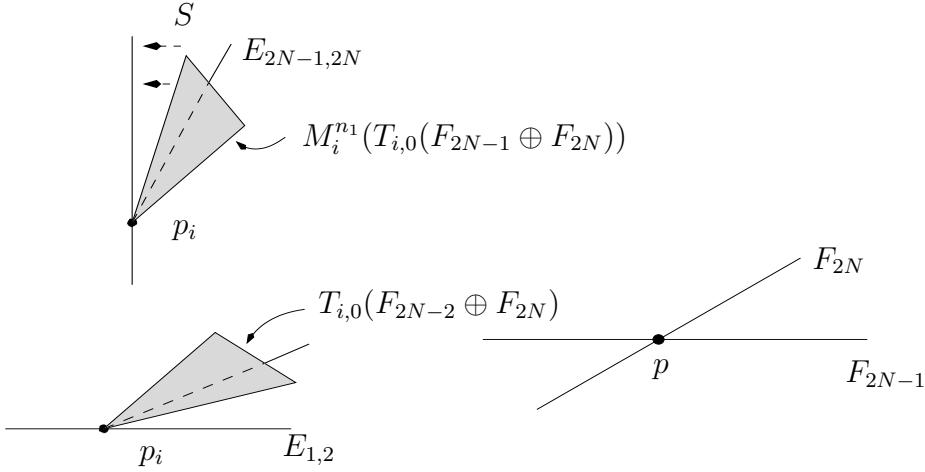}
  \caption{constructing a transition whose matrix sends
  $F_{2N-1} \oplus F_{2N}$ to  $E_{2N-1,2N}$ .}
  \label{fig1}
\end{figure}
By domination, for $n_1$ sufficiently large
$M_i^{n_1} \circ \, T_{i, 0} (F_{2N-1}\oplus F_{2N})$ is close
enough to $E_{2N-1, 2N}$.
Now by another small symplectic perturbation $S$ we may send
$M_i^{n_1} \circ T_{i, 0} (F_{2N-1} \oplus F_{2N})$ inside
$E_{2N-1, 2N}$.
Let $T_1 = S \circ M_i^{n_1} \circ T_{i, 0}$.
Now, using again the dominated splitting, there is $n_2$ sufficiently
large such that $ T_2^{-1}(E_{1,2}):= M^{-n_2} \circ T_1^{-1}(E_{1,2})$
is close enough to $F_1 \oplus F_2$.
It is possible to choose $v_1 , v_2 \in E_{1,2}$ such that
$$
T_2^{-1}(v_1) = f_1 + \sum_{j=3}^{2N} \varepsilon_j f_j  \quad
 \text{and} \quad
T_2^{-1}(v_2) = f_2 + \sum_{j=3}^{2N} \tilde{\varepsilon}_j f_j.
$$
Indeed, let $ V_{f_2}:= F_1 \oplus F_3 \cdots \oplus F_{2N}$.
As $\dim (V_{f_2}) = 2N-1$ and $\dim (T_1^{-1}(E_{1,2})) = 2$ there
exists $w:= f_1 + \varepsilon_3 f_3 + \cdots + \varepsilon_{2N} f_{2N}
\in V_{f_2} \cap T_1^{-1}(E_{1,2})$.
Since $w$ is close enough to $F_1 \oplus F_2$ it comes out that
$\varepsilon_i$, $i = 3, \cdots , 2N$ are arbitrarily small.
So, it is enough to take $v_1=  T_2(w)$.
Similarly we get $v_2$ satisfying the assertion.

Define a symplectic linear map $\tilde{I}$ as follows.
\begin{equation*}
\tilde{I}(f_i) =\left\{
\begin{array}{ll}
    T_2^{-1}(v_i)    & \text{if } i=1,2 \\
    f_{i}   & \text{if } i=2N, 2N-1 \\
    f_i + \varepsilon_{i^*}f_{2N} + \tilde{\varepsilon}_{i^*}f_{2N-1} &
      \text{if } 3 \leq i \leq N \\
    f_i - \varepsilon_{i^*} f_{2N} - \tilde{\varepsilon}_{i^*}f_{2N-1} &
      \text{if } N < i \leq 2N-2.
\end{array}
\right.
\end{equation*}
Let us check that $\tilde{I}$ is a symplectic transformation.
It is enough to verify that $\omega(\tilde{I}(f_j), \tilde{I}(f_k)) =
\omega(f_j, f_k)$ for all $1 \leq j, k \leq 2N, j \neq k$.
We have the following cases:

\smallskip

\noindent $(i)$ Let $j=1, 2$ and $3 \leq k \leq  N$.
One just prove for $j=1$ (the case  $j=2$ is analogous).
We have
\begin{align*}
\omega(\tilde{I}(f_1), \tilde{I}(f_{k})) & =
\omega (f_1 + \sum_{j=3}^{2N} \varepsilon_j f_j,
f_k + \varepsilon_{k^*}f_{2N} + \tilde{\varepsilon}_{k^*}f_{2N-1}) \\
& = \omega(f_1, \varepsilon_{k^*}f_{2N}) +
\omega (\varepsilon_{k^*} f_{k^*}, f_k) =
\varepsilon_{k^*} - \varepsilon_{k^*} \\
& = 0 = \omega(f_1,f_k).
\end{align*}

\noindent$(ii)$ Let $j = 1,2$ and $ N < k \le 2N-2$.
Again we just give the proof for $j=1$.
We have
\begin{align*}
\omega(\tilde{I}(f_1), \tilde{I}(f_{k})) & =
\omega(f_1 + \sum_{i=3}^{2N}\varepsilon_j f_j ,
f_k - \varepsilon_{k^*} f_{2N} - \tilde{\varepsilon}_{k^*}f_{2N-1}) \\
& = \omega(f_1, - \varepsilon_{k^*} f_{2N}) +
\omega(\varepsilon_{k^*} f_{k^*}, f_k) =
- \varepsilon_{k^*} + \varepsilon_{k^*} \\
& = 0 = \omega(f_1, f_k).
\end{align*}

\noindent $(iii)$ Let $3 \le j, k \le 2N-2$.
In this case we have
\begin{align*}
\omega(\tilde{I}(f_j) , \tilde{I}(f_k)) & =
\omega (f_j \pm \varepsilon_{j^*}f_{2N} \pm \tilde{\varepsilon}_{j^*}f_{2N-1} ,
f_k \pm \varepsilon_{k^*}f_{2N} \pm \tilde{\varepsilon}_{k^*}f_{2N-1}) \\
& = \omega(f_j, f_k).
\end{align*}

\noindent$(iv)$ Let $j = 1,2$ and $ k = 2N-1, 2N$.
Let us just show this case for the when $j=1$ and $k=2N$:
$$
\omega(\tilde{I}(f_1), \tilde{I}(f_{2N})) =
\omega (f_1 + \sum_{i=3}^{2N}\varepsilon_i f_i, f_{2N}) =
\omega(f_1, f_{2N}).
$$
This complete the proof that $\tilde{I}$ is symplectic.

Hence, $T:=T_2 \circ M_{p}^{n_2} \circ \tilde{I}$  maps $F_1 \oplus F_2$ and
$F_{2N-1} \oplus F_{2N}$ respectively to $E_{1,2}$ and  $E_{2N-1,2N}$.

It remains to prove that
$$
T (F_3 \oplus \cdots \oplus F_{2N-2}) = E_3 \oplus \cdots \oplus E_{2N-2}.
$$
Take $ 2 < j < 2N-1$ and assume that
$T (f_j)\notin E_3 \oplus \cdots \oplus E_{2N-2}$.
So, there is $v \in E_{1,2} \cup E_{2N-1, 2N}$ such that
$\omega( T(f_j) , v ) \neq 0$.
Let $w:= T^{-1}(v) \in (F_1 \oplus F_2) \cup (F_{2N-1} \oplus F_{2N})$.
Then,
$$
0= \omega(f_j,w) = \omega (T(f_j), T(w)) = \omega (T(f_j), v) \neq 0,
$$
which is a contradiction.
This completes the proof of the lemma.
\end{proof}

\begin{lemma}
\label{l.2}
Under the hypotheses of the proposition \ref{p.ultima}, if $i < 2N-1$
there exists a symplectic transition $T$ from $p$ to $p_i$ with the
following properties:
\begin{itemize}
\item  $T(F_{2N}) = E_{2N}$ and $T (F_{1} ) = E_{1},$
\item $T(F_2 \oplus \cdots \oplus F_{2N-1}) = E$, where
$E = E_2 \oplus \dots \oplus E_{2N-1}$ or
$E_{2,3} \oplus E_4 \oplus \dots \oplus E_{2N-3} \oplus E_{2N-2,2N-1}.$
\end{itemize}
\end{lemma}

\begin{proof}
Let $t_{i, 0}$ be an arbitrary symplectic $\varepsilon$-transition
from $p$ to $p_i$.
After a small symplectic perturbation we may suppose
$T_{i, 0} (F_{2N}) \notin E_1 \oplus E_2 \oplus \cdots \oplus E_{2N-1}$.
So, by the dominated splitting, for $n_1$ sufficiently large
$M_i^{n_1} \circ T_{i, 0} (F_{2N})$ is close enough to $E_{2N}$.
Now, by another small symplectic perturbation $S$ we may send
$M_i^{n_i} \circ T_{i, 0} (F_{2N})$ inside $E_{2N}$.
Let $T_1 = S \circ M_i^{n_1} \circ T_{i, 0}$.
Then $T_1(F_{2N}) = E_{2N}$.
All perturbations are symplectic, but we do not know whether
$T_1 (F_1) = E_1$ or not.

Let $\{f_1, \dots f_{2N}\}$ be a symplectic basis for $\mathcal{E}(P)$
and $e_{2N}=T_1(f_{2N}).$ For any $e_1 \in E_1$
we have
$$
\omega(T_1^{-1}(e_{2N}), T_1^{-1}(e_{1})) = \omega (e_{2N}, e_1)
\neq 0
$$
So, we conclude that $T_1^{-1}(E_{1})$ has a non null component in
the $F_1$ direction.
Now, using the dominated splitting, for $n_2$ sufficiently large
we have $ T_2^{-1}(e_1):= M^{-n_2} \circ T_1^{-1}(e_1)$ is close
enough to $f_1$ for some $e_1 \in E_1$.
More precisely for some $e_1 \in E_1$ we can write:
$$
T_2^{-1} (e_1) = f_1 + \varepsilon_2 f_2 + \varepsilon_3 f_3 +
\dots + \varepsilon_{2N} f_{2N}.
$$
where $f_j \in F_j$ and $\varepsilon_i$ are small enough whenever
$n_2$ is sufficiently large.

Now we define  a symplectic perturbation of the identity map
defined on the basis $\{f_1, f_2, \cdots , f_{2N}\}$ as follows.
\begin{equation*}
\tilde{I}(f_i) =\left\{
\begin{array}{ll}
    f_1 +
 \sum_{i=2}^{2N}\varepsilon_i f_i & \text{if } i=1 \\
    f_{2N} & \text{if } i=2N \\
    f_i + \varepsilon_{i^*}f_{2N}  & \text{if } 1< i \leq N \\
    f_i - \varepsilon_{i^*} f_{2N} & \text{if } N < i < 2N.
\end{array}
\right.
\end{equation*}

Let us verify that $\tilde{I}$ is a symplectic transformation.
It is enough to verify that
$\omega(\tilde{I}(f_j), \tilde{I}(f_k)) = \omega(f_j, f_k)$
for all $1 \leq j, k \leq 2N, j \neq k$.
We have the following cases:

\smallskip

\noindent $(i)$ Let $j=1 $ and $1 < k < N$.
Observe that for $k \leq N$ we have $\omega( f_{k^*}, f_k) = -1$.
So, we have
\begin{align*}
\omega(\tilde{I}(f_1), \tilde{I}(f_{k})) & =
  \omega(f_1 + \sum_{i=2}^{2N} \varepsilon_i f_i,f_{k} +
  \varepsilon_{k^*}f_{2N}) \\
& = \omega(f_1, \varepsilon_{k^*} f_{2N}) +
  \omega(\varepsilon_{k^*} f_{k^*}, f_k) \\
& = - \varepsilon_{k^*} + \varepsilon_{k^* }= 0 = \omega(f_1, f_k).
\end{align*}
$(ii)$ Let $j = 1, N < k < 2N$.
Note that for $k > N$, $\omega( f_{k^*}, f_k) = 1$.
Then,
\begin{align*}
\omega(\tilde{I}(f_1), \tilde{I}(f_{k}))& =
  \omega(f_1 +
  \sum_{i=2}^{2N}\varepsilon_i f_i,f_{k}-\varepsilon_{k^*}f_{2N}) \\
& = \omega(f_1, - \varepsilon_{k^*} f_{2N}) +
  \omega(\varepsilon_{k^*} f_{k^*},f_k) \\
& = \varepsilon_{k^*} - \varepsilon_{k^* }= 0 = \omega(f_1, f_k).
\end{align*}
$(iii)$ Let $1 < j, k < 2N$. In this case we have
$$
\omega(\tilde{I}(f_j) , \tilde{I}(f_k)) =
\omega ( f_j \pm \varepsilon_{j^*}f_{2N},
f_k \pm \varepsilon_{k^*}f_{2N}) = \omega(f_j, f_k).
$$
$(iv)$ Let $j = 1, k = 2N$. Then
$$
\omega(\tilde{I}(f_1), \tilde{I}(f_{2N})) =
\omega (f_1 + \sum_{i=2}^{2N}\varepsilon_i f_i, f_{2N}) =
\omega(f_1, f_{2N}).
$$

Hence, $T:= T_2 \circ \tilde{I}$ is a symplectic transition from $p$
to $p_i$ with $T (F_{2N}) = E_{2N}$ and $T(F_1) = E_1$.
Moreover, the fact that $T$ is symplectic implies that for all
$1<j<2N$
$$
\omega (T(f_j), T(f_{2N})) =  \omega (T(f_j), T(f_1)) = 0.
$$
This implies that $T(f_j)$ can not have coordinates in the $e_1$
and $e_{2N}$ directions.
Therefore, $T$ sends the eigenspace $F_2 \oplus \dots \oplus F_{2N-1}$
into $E$.

\end{proof}

\begin{proof}[Proof of Proposition \ref{p.ultima}.]
If $i=2N-1$ we apply Lemma \ref{l.1} and then use the method in
the proof of Lemma \ref{l.2} inductively to get the transition
$T_{i,0}$.

Otherwise, if $i < 2N-1$ we apply Lemma \ref{l.2} successively
$2N-i-1$ times to reduce the problem to the above case.

\end{proof}

\begin{proof} [Proof of Proposition \ref{p.identity}]
The proof of this proposition follows from the same arguments of
\cite[Lemma 5.4]{BDP03}.
Let us outline the main steps of the proof and point out how to
complete the proof in the symplectic case.

After a small symplectic perturbation we may suppose that $M_p$ is
diagonalizable for some $p \in \Sigma$ and there is $p_i \in \Sigma$
such that $p_i$ has complex eigenvalue of rank $(i, i+1)$.
Let
\begin{equation}
\mathcal{E}(p)= F_1 \prec \dots \prec F_{2N}
\end{equation}
and
\begin{equation} \label{decomp}
\mathcal{E}(p_i) = E_1 \prec \dots \prec E_{(i+1)^*, i^*} \prec
\dots \prec E_{i, i+1} \prec \dots \prec E_{2N}.
\end{equation}

By Proposition \ref{p.ultima} we can construct a symplectic
transition $[t^i]:= [t_{0,i}][t_{i,0}]$ from $p$ to itself whose matrix
$T^i$ preserves all subspaces $F_{j}$, $j \neq i, i+1, i^* , (i+1)^*$,
$F_i \oplus F_{i+1},$ and $F_{i^*} \oplus F_{(i+1)^*}$.

Using again Proposition \ref{p.ultima} and  the complex eigenvalues
inside $E_{i, i+1}$ we can obtain  a new transition $[\tilde{t}^i]$
such that $\tilde{T}^i$ interchanges $F_i$ and $F_{i+1}$.
This transition can be constructed exactly as in \cite[Lemma 5.7]{BDP03}
using Lemma \ref{pertu} and Lemma \ref{l.main}.
The unique difference is that to complete the proof in the symplectic
case we should prove that $\tilde{T}^i$ interchanges $F_{i^*}$ and
$F_{(i+1)^*}$ too.

Given a non-zero vector $f_{(i+1)} \in F_{(i+1)}$, let
$f_i = \tilde{T}^i (f_{(i+1)})$.
We show that $\tilde{T}^i (f_{i^*}) \in F_{(i+1)^*}$ for any
$f_{i^*} \in F_{i^*}.$

Since $F_{i^*} \oplus F_{(i+1)^*}$ is preserved by
$\tilde{T}^i$, there exists constants $r, s \in \mathbb{R}$
such that $ \tilde{T}^i(f_{i^*}) = r f_{i^*} + s f_{(i+1)^*}$.
Then,
\begin{align*}
0 & = \omega ( f_{i^*}, f_{(i+1)})
= \omega(\tilde{T}^i(f_{i^*}), \tilde{T}^i(f_{i+1})) \\
& = \omega (r f_{i^*} + s f_{(i+1)^*}, f_i)
= r \omega (f_{i^*}, f_i)
\end{align*}
The fact that $\omega (f_{i^*}, f_i) \neq 0$ implies $r=0$.
Therefore $\tilde{T}^i(f_{i^*})$ belongs to $F_{(i+1)^*}$.
Similarly we can prove that $\tilde{T}^i (F_{(i+1)^*}) = F_{i^*}$.

\end{proof}

\medskip

\noindent{\bf End of proof of Proposition \ref{p.bdp2.5}}.
In Proposition \ref{p.identity} we construct transitions
$[\tilde{t}^i]$ whose action on the finite set
$\{F_i(p)\}_{ 1 \leq j \leq 2N }$ of eigenspaces of $M_{\tilde{A}}(p)$
is  the transposition $(i, i+1)$ which interchanges $E_i(p)$ and
$E_{i+1}(p)$.
In what follows a combinatorial argument exactly as done in
\cite{BDP03} shows that after a small perturbation we obtain a
totally elliptic periodic point.
We verify that the arguments in their paper works also in the
symplectic case.

Given $0 \leq k < 2N$ denote by $\sigma_k$ the cyclic permutation
defined by $\sigma_k(E_j(p)) = E_{j+k}(p)$, where the sum $i+j$ is
considered in the cyclic group $\mathbb{Z}/(2N) \mathbb{Z}$.

As any permutation is a composition of transpositions, for every
$0 \le k < 2N$ there exists an element $[\tilde S_k]$ in the
semi-group generated by transitions $[\tilde t^i]$ such that if
$\tilde S_k$ is the matrix corresponding to the word
$[\tilde S_k]$ then one has $\tilde S_k (E_j(p)) = E_{j+k} (p)$.

Let $[S_k]$ be the word of matrices corresponding to the
perturbation $[\tilde S_k]$ in the semi-group generated by the
initials $[t^j]$.
As the $[t^j]$ are $\varepsilon_1$-transitions from $p$ to itself,
any word in the semi-group generated by the $[t^j]$, in particular
the $[S_k]$ is also an $\varepsilon_1$-transition from $p$ to
itself.
Let us write $[S_0] = [S_{2N}]$ the empty word whose corresponding
matrix is the identity.

By definition of transitions, for any $n \in \mathbb{N}$ there is
a point $x_n \in \Sigma$ such that the word $[M]_A(x_n)$ is
$\varepsilon_1$-close to the word $[W_n]$ corresponding to the
matrix $W_n$ defined by
$$
W_n :=  W_{2N-1, n} \circ \cdots W_{1, n} \circ  W_{0,n} ,
$$
where $W_{i, n} :=  S_{2N-i} \circ (M_{ A}(p))^n \circ
S_i \circ  S_{2N-i}  S_i$.

We know that for any $i$, the matrix $\tilde S_{2N-i} \circ \tilde
S_i$ acts trivially on the set of spaces $\{E_j(p)\}$.
Let us denote by
\begin{itemize}
\item  $\tilde \lambda_k$ the eigenvalues of $M_{\tilde A} (p)$
corresponding to  $E_k$,
\item $\mu_{i, j}$  the eigenvalue of $\tilde S_{2N-i} \circ \tilde S_i$
corresponding to the $E_j(p).$
\end{itemize}

Consequently, for every $j$ and any $n \in \mathbb{N}$ the space
$E_j(p)$ is an eigenspace of the matrix
$$
\tilde W_{i, n} := \tilde S_{2N-i} \circ (M_{\tilde A}(p))^n \circ
\tilde S_i \circ \tilde S_{2N-i} \tilde S_i,
$$
whose corresponding eigenvalue is $\mu_{i,j}^2 \tilde \lambda_{i+j}^n$.

As the transitions are symplectic we have
$\tilde \lambda_i \tilde \lambda_{i^*}= 1$ and
$\mu_{i,j} \mu_{i,j^*} = 1$.
It comes out that
$$
\prod_{i=1}^{2N} \tilde \lambda_i = 1
$$
and
$$
 C_j:=\prod_{i=0}^{2N-1} \mu_{i,j}^2 = (\prod_{i=0}^{2N-1}
\mu_{i,j^*}^2)^{-1}= C_{j^*}^{-1}.
$$

The word $[\tilde{W}_n]$ corresponding to the matrix defined as
$$
\tilde{W}_n :=  \tilde{W}_{2N-1, n} \circ \cdots  \tilde{W}_{1, n}
\circ  \tilde{W}_{0,n}
$$
is $(\varepsilon_0 + \varepsilon_1)$-close to $[W_n]$ and so it
is an $\varepsilon_0 + 2 \varepsilon_1$-perturbation of the word
$[M]_A(x_n)$.
So we conclude that $E_j(p)$ is an eigenspace of $\tilde W_n$ with
eigenvalue $C_j$.

Observe that $C_j$ are not necessarily close to $1$.
Consider $B_n$ matrices having $E_j(p)$ as eigenspaces and
$(C_j)^{-\frac{1}{n}}$ as their eigenvalues.
Observe that $B_n$ is symplectic too and $(B_n)^n = \tilde{W}_n^{-1}$.
Denote by $[M]_{\hat{A}}(p)$ the word obtained from
$[M]_{\tilde {A}}(p)$ by replacing its first letter $\tilde A(p)$
by $\tilde A(p) \circ B_n $.
For $n$ large enough this new word is an $\epsilon_1$-perturbation of
$[M]_{\tilde {A}}(p)$, so by item $(i)$ of Proposition \ref{p.identity}
it is also $2\epsilon_1$-perturbation of $[M]_{A}(p)$.
As $B_n$ commutes with $M_{\tilde {A}}(p)$
we get
$$
(M_{\tilde {A}}(p) \circ B_n)^n = M_{\tilde {A}}^n(p) \circ
\tilde{W}_n^{-1}.
$$

So, the word $[\hat W_n]$ obtained by changing the initial subword
$[M]_{\tilde {A}}^n(p)$ of $[\tilde W_n]$ by
$[M]_{\hat {A}}(p)$ is $(\epsilon_0 + 2 \epsilon_1) < \epsilon$ close to
the word $[M]_{A}(x_n)$ and its corresponding matrix
$\hat W_n = \tilde W_n \circ \tilde{W}_n^{-1} = Id$.
This completes the proof.

\bibliographystyle{plain}

\bigskip

\flushleft

{\bf Vanderlei Horita} (vhorita\@@mat.ibilce.unesp.br)\\
UNESP - São Paulo State Universtiy \\
Departamento de Matem\'atica, IBILCE \\
Rua Crist\'ov\~ao Colombo 2265 \\
15054-000 S. J. Rio Preto, SP, Brazil

\bigskip

\flushleft

{\bf Ali Tahzibi} (tahzibi\@@icmc.usp.br) \\
Departamento de Matemática e de Computação, ICMC/USP\\
Av. Trabalhador São Carlense , 400-Cx. Postal 668 \\
13560-320 São Carlos, SP, Brazil\\
 www.icmc.usp.br/$\sim$tahzibi
\end{document}